\documentclass[11pt]{article}

\usepackage{latexsym}
\usepackage{amssymb}

\newtheorem{thm}{Theorem}

\begin{document}
{
\begin{center}
{\Large\bf
The truncated matrix trigonometric moment problem: the operator approach.}
\end{center}
\begin{center}
{\bf S.M. Zagorodnyuk}
\end{center}

\section{Introduction.}
The aim of this investigation is to obtain a bijective parameterization of all solutions of the
truncated matrix trigonometric moment problem.
Recall that the truncated matrix trigonometric moment problem
consists of finding a non-decreasing $\mathbb{C}_{N\times N}$-valued function
$M(t) = (m_{k,l})_{k,l=0}^{N-1}$, $t\in [0,2\pi]$, $M(0)=0$, which is
left-continuous in $(0,2\pi]$, and such that
\begin{equation}
\label{f1_1}
\int_0^{2\pi} e^{int} dM(t) = S_n,\qquad n=0,1,...,d,
\end{equation}
where $\{ S_n \}_{n=0}^d$ is a prescribed sequence of $(N\times N)$ complex matrices (moments).
Here $N\in \mathbb{N}$ and $d\in \mathbb{Z}_+$ are fixed numbers.

\noindent
Set
\begin{equation}
\label{f1_2}
T_d = (S_{i-j})_{i,j=0}^d =
\left(
\begin{array}{ccccc} S_0 & S_{-1} & S_{-2} & \ldots & S_{-d}\\
S_1 & S_0 & S_{-1} & \ldots & S_{-d+1}\\
S_2 & S_1 & S_0 & \ldots & S_{-d+2}\\
\vdots & \vdots & \vdots & \ddots & \vdots\\
S_d & S_{d-1} & S_{d-2} & \ldots & S_0\end{array}
\right),
\end{equation}
where
\begin{equation}
\label{f1_3}
S_n := S_{-n}^*,\qquad n=-d,-d+1,...,-1.
\end{equation}
The scalar ($N=1$) truncated trigonometric moment problem is well investigated.
In~1911, Riesz and Herglotz  obtained the necessary and sufficient conditions of the solvability for this moment
problem (e.g.~\cite{cit_900_Akh}).
Canonical solutions of the moment problem were described in the positive definite case:
$T_d>0$, by Krein and~Nudelman~\cite{cit_1000_KN}.
In~1966, Chumakin described all solutions of the scalar truncated trigonometric moment
problem using his results on the generalized resolvents of isometric operators,
see~\cite{cit_2000_Ch},\cite{cit_3000_Ch},\cite{cit_4000_Ch}.

In the general case of  an arbitrary $N$, the following condition:
\begin{equation}
\label{f1_4}
T_d\geq 0,
\end{equation}
is necessary and sufficient for the solvability of the moment problem~(\ref{f1_1}) (e.g.~\cite{cit_4500_A}).
In~1969, Inin obtained a description of all solutions of the truncated matrix trigonometric moment
problem in the completely indeterminate case: $T_d>0$~\cite{cit_5000_In}. He used the tools of the theory of
pseudo-Hilbert spaces developed by Krein and Berezanskii~\cite{cit_6000_M}.

The truncated matrix trigonometric moment problem is closely related to (and it is essentially the same as)
the Carath\'eodory matrix coefficient problem (this relation is based on the matrix
extension of the Riesz-Herglotz integral representation), see~\cite{cit_6500_ChH},\cite{cit_6700_FK} and
References therein.
In~1998, a parameterization of
all solutions of the last problem both in nondegenerate and degenerate cases was for the first time obtained
by Chen and Hu~\cite{cit_6500_ChH}. In~2006, another parameterization of
all solutions of this problem both in nondegenerate and degenerate cases was obtained
by Fritzsche and Kirstein~\cite{cit_6700_FK}.
However, it is not clear whether the above parameterizations are bijective.

\noindent
We shall describe all solutions of the truncated matrix trigonometric moment
problem in the general case: $T_d\geq 0$, as well.
We shall use an abstract operator approach and the mentioned above results of Chumakin on
generalized resolvents of isometric operators. The abstract operator approach allows investigate simultaneously
both the nondegenerate and degenerate cases of different moment problems~(see~\cite{cit_7000_Z} and
\cite{cit_8000_Z},\cite{cit_9000_Z}).
The obtained parameterization of all solutions is bijective. Also, the operator point of view
allows to see the transparent whole picture of the problem at once rather then step-by-step algorithm.


{\bf Notations. }
As usual, we denote by $\mathbb{R}, \mathbb{C}, \mathbb{N}, \mathbb{Z}, \mathbb{Z}_+$,
the sets of real numbers, complex numbers, positive integers, integers and non-negative integers,
respectively; $\mathbb{D} = \{ z\in \mathbb{C}:\ |z|<1 \}$.
The set of all complex vectors of size $N$: $a = (a_0,a_1,\ldots,a_{N-1})$, we denote by
$\mathbb{C}_N$, $N\in \mathbb{N}$.
If $a\in \mathbb{C}^N$, then $a^*$ means its complex conjugate vector.
The set of all complex matrices of size $(N\times N)$ we denote by $\mathbb{C}_{N\times N}$.

\noindent
Let $M(x)$ be a left-continuous non-decreasing matrix function $M(x) = ( m_{k,l}(x) )_{k,l=0}^{N-1}$
on $[0,2\pi]$, $M(0)=0$, and $\tau_M (x) := \sum_{k=0}^{N-1} m_{k,k} (x)$;
$\Psi(x) = ( dm_{k,l}/ d\tau_M )_{k,l=0}^{N-1}$.  By $L^2(M)$ we denote a set (of classes of equivalence)
of $\mathbb{C}_N$-valued functions $f$ on $[0,2\pi]$,
$f = (f_0,f_1,\ldots,f_{N-1})$, such that (see, e.g.,~\cite{cit_10000_MM})
$$ \| f \|^2_{L^2(M)} := \int_0^{2\pi}  f(x) \Psi(x) f^*(x) d\tau_M (x) < \infty. $$
The space $L^2(M)$ is a Hilbert space with a scalar product
$$ ( f,g )_{L^2(M)} := \int_0^{2\pi}  f(x) \Psi(x) g^*(x) d\tau_M (x),\qquad f,g\in L^2(M). $$

If H is a Hilbert space then $(\cdot,\cdot)_H$ and $\| \cdot \|_H$ mean
the scalar product and the norm in $H$, respectively.
Indices may be omitted in obvious cases.

\noindent For a linear operator $A$ in $H$, we denote by $D(A)$
its  domain, by $R(A)$ its range, by $\mathop{\rm Ker}\nolimits A$
its null subspace (kernel), and $A^*$ means the adjoint operator
if it exists. If $A$ is invertible then $A^{-1}$ means its
inverse. $\overline{A}$ means the closure of the operator, if the
operator is closable. If $A$ is bounded then $\| A \|$ denotes its
norm.
For a set $M\subseteq H$
we denote by $\overline{M}$ the closure of $M$ in the norm of $H$.
For an arbitrary set of elements $\{ x_n \}_{n\in I}$ in
$H$, we denote by $\mathop{\rm Lin}\nolimits\{ x_n \}_{n\in I}$
the set of all linear combinations of elements $x_n$,
and $\mathop{\rm span}\nolimits\{ x_n \}_{n\in I}
:= \overline{ \mathop{\rm Lin}\nolimits\{ x_n \}_{n\in I} }$.
Here $I$ is an arbitrary set of indices.
By $E_H$ we denote the identity operator in $H$, i.e. $E_H x = x$,
$x\in H$. If $H_1$ is a subspace of $H$, then $P_{H_1} =
P_{H_1}^{H}$ is an operator of the orthogonal projection on $H_1$
in $H$.

\section{Solvability and a description of solutions.}
Let the moment problem~(\ref{f1_1}) be given. Suppose that the moment problem has a solution $M(t)$.
Define $S_n$, $n=-d,-d+1,...,-1$, by relation~(\ref{f1_3}).
Observe that
$$ S_{-n} = S_n^* = \left( \int_0^{2\pi} e^{int} dM(t) \right)^* =
\int_0^{2\pi} e^{-int} dM(t),\qquad n=0,1,...,d. $$
Therefore
\begin{equation}
\label{f2_1}
\int_0^{2\pi} e^{int} dM(t) = S_n,\qquad -d\leq n\leq d.
\end{equation}
Consider an arbitrary vector-valued polynomial $P(t)$ of the following form:
\begin{equation}
\label{f2_2}
P(t) = \sum_{k=0}^d (\alpha_{k,0},\alpha_{k,1},...,\alpha_{k,N-1}) e^{ikt} =
\sum_{k=0}^d \sum_{s=0}^{N-1} \alpha_{k,s} e^{ikt} \vec{e}_s,\ \alpha_{k,s}\in \mathbb{C},
\end{equation}
where
\begin{equation}
\label{f2_3}
\vec{e}_s = (\delta_{0,s},\delta_{1,s},...,\delta_{N-1,s}).
\end{equation}
Then
$$ 0\leq \int_0^{2\pi} P(t) dM(t) P^*(t) =
\sum_{k,r=0}^d \sum_{s,l=0}^{N-1} \alpha_{k,s} \overline{\alpha_{r,l}} \int_0^{2\pi} e^{i(k-r)t} \vec{e}_s
dM(t) \vec{e}_l^* $$
$$ =
\sum_{k,r=0}^d \sum_{s,l=0}^{N-1} \alpha_{k,s} \overline{\alpha_{r,l}} \vec{e}_s S_{k-r}
\vec{e}_l^*.
$$
Thus, for arbitrary complex numbers $\alpha_{k,s}$, $0\leq k\leq d$, $0\leq s\leq N-1$, it holds
\begin{equation}
\label{f2_4}
\sum_{k,r=0}^d \sum_{s,l=0}^{N-1} \alpha_{k,s} \overline{\alpha_{r,l}} \vec{e}_s S_{k-r}
\vec{e}_l^* \geq 0.
\end{equation}
Define the matrix $T_d$ by~(\ref{f1_2}) and set
\begin{equation}
\label{f2_5}
\vec\alpha =
(\alpha_{0,1},\alpha_{0,2},...,\alpha_{0,N},\alpha_{1,1},\alpha_{1,2},...,\alpha_{1,N},
...,\alpha_{d,1},\alpha_{d,2},...,\alpha_{d,N}).
\end{equation}
By the rules of the multiplication of block matrices relation~(\ref{f2_4}) means that
\begin{equation}
\label{f2_6}
\vec\alpha T_d \vec\alpha^* \geq 0,
\end{equation}
and therefore
\begin{equation}
\label{f2_7}
T_d\geq 0.
\end{equation}
Conversely, let the moment problem~(\ref{f1_1}) with $d\in \mathbb{N}$ be given and relation~(\ref{f2_7}) holds
with $T_d$ defined by~(\ref{f1_2}).
Let
\begin{equation}
\label{f2_8}
T_d = (\gamma_{n,m})_{n,m=0}^{(d+1)N-1},\
S_k = ( S_{k;s,l} )_{s,l=0}^{N-1},\quad -d\leq k\leq d,
\end{equation}
where $\gamma_{n,m}, S_{k;s,l}\in \mathbb{C}$.

\noindent
Observe that
\begin{equation}
\label{f2_9}
\gamma_{kN+s,rN+l} = S_{k-r;s,l},\qquad 0\leq k,r\leq d,\quad 0\leq s,l\leq N-1.
\end{equation}
By the well-known construction~(see, e.g.,~\cite[Theorem 1]{cit_9000_Z}), from~(\ref{f2_7}) it follows that
there exist a Hilbert space $H$ and elements $\{ x_n \}_{n=0}^{ (d+1)N-1 }$ in $H$ such that
\begin{equation}
\label{f2_10}
(x_n,x_m) = \gamma_{n,m},\qquad 0\leq n,m\leq (d+1)N-1,
\end{equation}
and $\mathop{\rm span}\nolimits\{ x_n \}_{n=0}^{ (d+1)N-1 } = H$.
Set $H_0 := \mathop{\rm Lin}\nolimits\{ x_n \}_{n=0}^{ dN-1 }$.

\noindent
Consider the following operator:
\begin{equation}
\label{f2_11}
A x = \sum_{k=0}^{dN-1} \alpha_k x_{k+N},\quad x = \sum_{k=0}^{dN-1} \alpha_k x_k,\ \alpha_k\in \mathbb{C}.
\end{equation}
Let us check that this definition is correct. Suppose that $x\in H_0$ has two representations:
$$ x = \sum_{k=0}^{dN-1} \alpha_k x_k,\quad x = \sum_{k=0}^{dN-1} \beta_k x_k,\quad \alpha_k,\beta_k\in
\mathbb{C}. $$
Using relations~(\ref{f2_9}),(\ref{f2_10}) we may write
$$ \left\| \sum_{k=0}^{dN-1} \alpha_k x_{k+N} - \sum_{k=0}^{dN-1} \beta_k x_{k+N} \right\|^2 =
\left(
\sum_{k=0}^{dN-1} (\alpha_k-\beta_k) x_{k+N}, \sum_{r=0}^{dN-1} (\alpha_r-\beta_r) x_{r+N}
\right) $$
$$ = \sum_{k,r=0}^{dN-1} (\alpha_k-\beta_k) \overline{ (\alpha_r-\beta_r) } \gamma_{k+N,r+N}
= \sum_{k,r=0}^{dN-1} (\alpha_k-\beta_k) \overline{ (\alpha_r-\beta_r) } \gamma_{k,r} $$
$$ =
\left(
\sum_{k=0}^{dN-1} (\alpha_k-\beta_k) x_{k}, \sum_{r=0}^{dN-1} (\alpha_r-\beta_r) x_{r}
\right) = (x-x,x-x) = 0. $$
Thus, the definition of $A$ is correct and  $D(A) = H_0$.

\noindent
Let $x,y\in H_0$, $x = \sum_{k=0}^{dN-1} \alpha_k x_k$, $y =\sum_{k=0}^{dN-1} \gamma_k x_k$,
$\alpha_k,\gamma_k\in \mathbb{C}$.
Then
$$ (A x, A y) = \sum_{k,r=0}^{dN-1} \alpha_k \overline{ \gamma_r } (x_{k+N},x_{r+N}) =
\sum_{k,r=0}^{dN-1} \alpha_k \overline{ \gamma_r } \gamma_{k+N,r+N} $$
$$ = \sum_{k,r=0}^{dN-1} \alpha_k \overline{ \gamma_r } \gamma_{k,r}
= \sum_{k,r=0}^{dN-1} \alpha_k \overline{ \gamma_r } (x_k,x_r) = (x,y). $$
Thus, $A$ is an isometric operator. Every isometric operator admits a unitary extension (\cite{cit_11000_AG}).
Let $U\supseteq A$ be a unitary extension of $A$ in a Hilbert space $\widetilde H\supseteq H$.
Choose an arbitrary non-negative integer $n$:
$$ n = rN+l,\qquad  0\leq r\leq d,\ 0\leq l\leq N-1. $$
By induction one easily derives the following  relation:
\begin{equation}
\label{f2_12}
x_{rN+l} = A^r x_l.
\end{equation}
Choose an arbitrary $m$:
$$ m = kN+s,\qquad  0\leq k\leq d,\ 0\leq s\leq N-1. $$
Using~(\ref{f2_9}) we may write
$$ S_{k-r;s,l} = \gamma_{kN+s,rN+l} = (x_m,x_n)_H = (A^k  x_s, A^r x_l)_H =
(U^k  x_s, U^r x_l)_{\widetilde H} $$
$$ = (U^{k-r}  x_s, x_l)_{\widetilde H} =
\int_0^{2\pi} e^{ i(k-r)t } d (E_t x_s,x_l)_{\widetilde H}, $$
where $\{ E_t \}_{t\in [0,2\pi]}$ is the left-continuous orthogonal resolution of unity of the operator $U$.
Thus, we have
\begin{equation}
\label{f2_13}
S_{j;s,l} = \int_0^{2\pi} e^{ ijt } d (P^{\widetilde H}_H E_t x_s,x_l)_H,\quad -d\leq j\leq d,\ 0\leq s,l\leq N-1.
\end{equation}
Set
\begin{equation}
\label{f2_14}
M_U (t) = \left( (P^{\widetilde H}_H E_t x_s,x_l)_H \right)_{s,l=0}^{N-1},\qquad t\in [0,2\pi].
\end{equation}
Then $M_U(t)$ is a solution of the moment problem~(\ref{f1_1}) (the fact that it is non-decreasing
follows easily from the properties of the orthogonal resolution of unity).

We can state the following well-known fact:
\begin{thm}
\label{t2_1}
Let the truncated matrix trigonometric moment problem~(\ref{f1_1}) be given.
The moment problem has a solution if and only if relation~(\ref{f1_4}) with
$T_d$ from~(\ref{f1_2}) holds.
\end{thm}
{\bf Proof. } The required result for the case $d\in \mathbb{N}$ was proved above.
For the case $d=0$ the following function is a solution:
$$ M(t) = \left\{
\begin{array}{cc}
0, & t=0 \\
S_0, & t\in (0,2\pi]\end{array}
\right.. $$
$\Box$

{\bf Remark. }
For the case $d=0$, an arbitrary
left-continuous non-decreasing $\mathbb{C}_{N\times N}$-valued function $M$,
$M(0)=0$, $M(2\pi) = S_0$ is a solution of the moment problem~(\ref{f1_1}). Therefore we
shall investigate the case $d\in \mathbb{N}$.

Let the moment problem~(\ref{f1_1}) be given with $d\in \mathbb{N}$ and condition~(\ref{f1_4}) holds.
As it was done above, we construct a Hilbert space $H$, a sequence $\{ x_n \}_{n=0}^{ (d+1)N-1 }$ in $H$ and
the isometric operator $A$. Let $\widehat U\supseteq A$ be an arbitrary unitary extension of
$A$ in a Hilbert space $\widehat H\supseteq H$. Let $\{ \widehat E_t \}_{t\in [0,2\pi]}$ be
the left-continuous orthogonal resolution of unity of $\widehat U$. Recall~\cite{cit_3000_Ch},\cite{cit_4000_Ch}
that the following function:
\begin{equation}
\label{f2_15}
\mathbf{E}_t = P^{\widehat H}_H \widehat E_t,\qquad t\in [0,2\pi],
\end{equation}
is said to be a {\it spectral function} of $A$. The operator-valued function
\begin{equation}
\label{f2_16}
\mathbf{R}_\zeta = P^{\widehat H}_H (E_{\widehat H} - \zeta \widehat U)^{-1},\qquad \zeta\in \mathbb{C}:
|\zeta|\not= 1,
\end{equation}
is said to be a {\it generalized resolvent} of $A$.
If $\mathbf{E}_t$ and $\mathbf{R}_\zeta$ correspond to the same unitary extension of $A$, they
are said to be related. The related left-continuous spectral function and generalized resolvent of $A$ are
in a bijective correspondence:
\begin{equation}
\label{f2_17}
(\mathbf{R}_\zeta h,g)_H = \int_0^{2\pi} \frac{1}{1-\zeta e^{it}} d(\mathbf{E}_t h,g),
\qquad \forall h,g\in H.
\end{equation}
The function $(\mathbf{E}_t h,g)$  can be found by the
inversion formula (\cite{cit_4000_Ch}).

As we have seen above, an arbitrary left-continuous spectral function of the isometric operator $A$ generates
a solution of the moment problem~(\ref{f1_1}) by relation~(\ref{f2_14}).

On the other hand, let $\widehat M$ be an arbitrary solution of the moment problem~(\ref{f1_1}).
A set of all (classes of equivalence which include) polynomials of the form~(\ref{f2_2})
in $L^2(\widehat M)$ we shall denote by $L^2_{0,d}(\widehat M)$.
Choose an arbitrary
\begin{equation}
\label{f2_18}
Q(t) =
\sum_{r=0}^d \sum_{l=0}^{N-1} \beta_{r,l} e^{irt} \vec{e}_l,\ \beta_{r,l}\in \mathbb{C}.
\end{equation}
Then
$$ (P(t),Q(t))_{ L^2(\widehat M) } =
\sum_{k,r=0}^d \sum_{s,l=0}^{N-1} \alpha_{k,s} \overline{\beta_{r,l}}
\int_0^{2\pi} e^{ i(k-r)t } \vec{e}_s d\widehat M(t) \vec{e}_l^* $$
$$ = \sum_{k,r=0}^d \sum_{s,l=0}^{N-1} \alpha_{k,s} \overline{\beta_{r,l}}
 \vec{e}_s S_{k-r} \vec{e}_l^*
= \sum_{k,r=0}^d \sum_{s,l=0}^{N-1} \alpha_{k,s} \overline{\beta_{r,l}}
 S_{k-r;s,l} $$
$$ = \sum_{k,r=0}^d \sum_{s,l=0}^{N-1} \alpha_{k,s} \overline{\beta_{r,l}}
 \gamma_{kN+s,rN+l}
= \sum_{k,r=0}^d \sum_{s,l=0}^{N-1} \alpha_{k,s} \overline{\beta_{r,l}}
 (x_{kN+s}, x_{rN+l})_H $$
 \begin{equation}
\label{f2_19}
= \left(
\sum_{k=0}^d \sum_{s=0}^{N-1} \alpha_{k,s} x_{kN+s},
\sum_{r=0}^d \sum_{l=0}^{N-1} \beta_{r,l} x_{rN+l}
\right)_H.
\end{equation}
Consider the following operator:
\begin{equation}
\label{f2_20}
W P(t) = \sum_{k=0}^d \sum_{s=0}^{N-1} \alpha_{k,s} x_{kN+s}.
\end{equation}
Let us check that this operator is correctly defined as an operator from
$L^2_{0,d}(\widehat M)$ to $H$. Let $P(t)$ and $Q(t)$ are two polynomials
of the forms~(\ref{f2_2}) and~(\ref{f2_18}), respectively.
Suppose that they belong to the same class of equivalence in $L^2(\widehat M)$:
\begin{equation}
\label{f2_21}
(P(t)-Q(t),P(t)-Q(t))_{L^2(\widehat M)} = 0.
\end{equation}
Then
$$ 0 =
\left(
\sum_{k=0}^d \sum_{s=0}^{N-1} (\alpha_{k,s}-\beta_{k,s}) e^{ikt} \vec e_s,
\sum_{r=0}^d \sum_{l=0}^{N-1} (\alpha_{r,l}-\beta_{r,l}) e^{irt} \vec e_l
\right)_{L^2(\widehat M)} $$
$$ =
\sum_{k,r=0}^d \sum_{s,l=0}^{N-1} (\alpha_{k,s}-\beta_{k,s}) \overline{ (\alpha_{r,l}-\beta_{r,l}) }
\int_0^{2\pi} e^{ i(k-r)t } \vec{e}_s d\widehat M(t) \vec{e}_l^* $$
$$ = \left(
\sum_{k=0}^d \sum_{s=0}^{N-1} (\alpha_{k,s}-\beta_{k,s}) x_{kN+s},
\sum_{r=0}^d \sum_{l=0}^{N-1} (\alpha_{r,l}-\beta_{r,l}) x_{rN+l}
\right)_H = \| WP - WQ \|_H. $$
Thus, the operator $W$ is defined correctly.
Relation~(\ref{f2_19}) shows that $W$ is an isometric operator. It maps
$L^2_{0,d}(\widehat M)$ onto $H$.
Denote
$$ L^2_1(\widehat M) := L^2(\widehat M) \ominus L^2_{0,d}(\widehat M). $$
The operator
$$ U := W \oplus E_{ L^2_1(\widehat M) }, $$
is a unitary operator which maps $L^2(\widehat M) = L^2_{0,d}(\widehat M) \oplus L^2_1(\widehat M)$
onto $H_1 := H\oplus L^2_1(\widehat M)$.

\noindent
Consider the following unitary operator:
$$ U_0 f(t) = e^{it} f(t),\qquad f(t) \in L^2(\widehat M). $$
Then
$$ \widetilde U_0 := U U_0 U^{-1}, $$
is a unitary operator in $H_1$. Observe that
$$ \widetilde U_0 x_{kN+s} = U U_0 e^{ikt} \vec{e}_s = U e^{i(k+1)t} \vec{e}_s =
x_{(k+1)N+s} = A x_{kN+s}, $$
where $0\leq k\leq d-1,\ 0\leq s\leq N-1$.
Therefore $\widetilde U_0\supseteq A$.
Let $\{ \widetilde E_t \}_{ t\in [0,2\pi] }$ be the left-continuous orthogonal
resolution of unity of $\widetilde U_0$,
and $\mathbf{E}_t$, $\mathbf{R}_\zeta$, be a spectral function and a generalized resolvent of $A$
which correspond to the unitary extension $\widetilde U_0$, respectively.
Let us check that
\begin{equation}
\label{f2_22}
\widehat M(t) = ( (\mathbf{E}_t x_s,x_l)_H )_{s,l=0}^{N-1}.
\end{equation}
In fact, we may write
$$ \int_0^{2\pi} \frac{1}{1-\zeta e^{it}} d(\mathbf{E}_t x_s,x_l)_H =
(\mathbf{R}_\zeta x_s,x_l)_H =
\left(
( E_{H_1} - \zeta \widetilde U_0 )^{-1} x_s, x_l
\right)_{H_1} $$
$$ =
\left(
U ( E_{L^2(\widehat M)} - \zeta U_0 )^{-1} U^{-1} x_s, x_l
\right)_{H_1}
=
\left(
( E_{L^2(\widehat M)} - \zeta U_0 )^{-1} \vec e_s, \vec e_l
\right)_{ L^2(\widehat M) } $$
$$ =
\int_0^{2\pi} \frac{1}{1-\zeta e^{it}} \vec e_s d\widehat M(t) \vec e_l^*. $$
By the inversion formula we conclude that relation~(\ref{f2_22}) is true.

\begin{thm}
\label{t2_2}
Let the truncated matrix trigonometric moment problem~(\ref{f1_1}) with $d\in \mathbb{N}$ be given and
condition~(\ref{f1_4}) is true. Let an operator $A$ be constructed for the
moment problem as in~(\ref{f2_11}).
All solutions of the moment problem have the following form
\begin{equation}
\label{f2_23}
M(t) = (m_{k,j} (t))_{k,j=0}^{N-1},\quad
m_{k,j} (t) = ( \mathbf E_t x_k, x_j)_H,
\end{equation}
where $\mathbf E_t$ is a left-continuous spectral function of the isometric operator $A$.
Conversely, an arbitrary left-continuous spectral function of $A$ generates by formula~(\ref{f2_23})
a solution of the moment problem~(\ref{f1_1}).

Moreover, the correspondence between all left-continuous spectral functions of $A$ and all solutions
of the moment problem is bijective.
\end{thm}
{\bf Proof. }
It remains to prove that different left-continuous spectral functions of the operator $A$ produce different
solutions of the moment problem~(\ref{f1_1}).
Set
\begin{equation}
\label{f2_24}
H_\zeta := (E_H - \zeta A) D(A) = (E_H - \zeta A) H_0,\quad \zeta\in \mathbb{C}: | \zeta |\not= 1;
\end{equation}
\begin{equation}
\label{f2_25}
L_N := \mathop{\rm Lin}\nolimits\{ x_k \}_{k=0}^{N-1}.
\end{equation}
Choose an arbitrary element $x\in H$, $x=\sum_{k=0}^{dN+N-1} \alpha_k x_k$, $\alpha_k\in \mathbb{C}$.
Let us check that for an arbitrary $\zeta\in \mathbb{C}\backslash\{ 0\}: | \zeta |\not= 1$,
there exists a representation
\begin{equation}
\label{f2_26}
x = v + y,\qquad v\in H_\zeta,\ y\in L_N,
\end{equation}
where elements $v,y$ may depend on the choice of $\zeta$.

\noindent
In fact, choose an arbitrary $\zeta\in \mathbb{C}\backslash\{ 0\}: |z|\not= 1$. Set
\begin{equation}
\label{f2_27}
c_r := -\frac{1}{\zeta} \alpha_{r+N},\qquad r = dN-N, dN-N+1,...,dN-1.
\end{equation}
Then we set
\begin{equation}
\label{f2_28}
c_r := \frac{1}{\zeta} (c_{r+N} - \alpha_{r+N}),\qquad r = dN-N-1, dN-N-2,...,0.
\end{equation}
Let
\begin{equation}
\label{f2_29}
u := \sum_{k=0}^{dN-1} c_k x_k \in D(A);
\end{equation}
\begin{equation}
\label{f2_30}
v := (E_H - \zeta A) u \in H_\zeta.
\end{equation}
Then
$$ v = \sum_{k=0}^{dN-1} c_k x_k - \zeta \sum_{k=0}^{dN-1} c_k x_{k+N}
 = \sum_{k=0}^{dN-1} c_k x_k - \zeta \sum_{k=N}^{dN+N-1} c_{k-N} x_k $$
$$ = \sum_{k=0}^{N-1} c_k x_k + \sum_{k=N}^{dN-1} ( c_k - \zeta c_{k-N} ) x_k
- \zeta \sum_{k=dN}^{dN+N-1} c_{k-N} x_k $$
$$ = \sum_{k=0}^{N-1} c_k x_k + \sum_{k=N}^{dN+N-1} \alpha_k x_k
 = \sum_{k=0}^{N-1} (c_k-\alpha_k) x_k + x. $$
Finally, we set $y := - \sum_{k=0}^{N-1} (c_k-\alpha_k) x_k \in L_N$, and obtain
$x = v + y$. Thus, relation~(\ref{f2_26}) holds.

Suppose to the contrary that two different left-continuous spectral functions of $A$ produce the same solution of
the moment problem~(\ref{f1_1}). That means that
there exist two unitary extensions
$U_j\supseteq A$, in Hilbert spaces $\widetilde H_j\supseteq H$, such that
\begin{equation}
\label{f2_30_1}
\mathbf{E}_{1,t} = P_{H}^{\widetilde H_1} E_{1,t} \not= P_{H}^{\widetilde H_2} E_{2,t} = \mathbf{E}_{2,t};
\end{equation}
and
\begin{equation}
\label{f2_31}
(P_{H}^{\widetilde H_1} E_{1,t} x_k,x_j)_H = (P_{H}^{\widetilde H_2} E_{2,t} x_k,x_j)_H,\qquad
0\leq k,j\leq N-1,\quad t\in [0,2\pi],
\end{equation}
where $\{ E_{j,t} \}_{t\in [0,2\pi]}$ are orthogonal resolutions of unity of
operators $U_j$, $j=1,2$.
By linearity we get
\begin{equation}
\label{f2_32}
(P_{H}^{\widetilde H_1} E_{1,t} x,y)_H = (P_{H}^{\widetilde H_2} E_{2,t} x,y)_H,\qquad
x,y\in L_N,\quad t\in [0,2\pi].
\end{equation}
Set
$$ R_{j,\zeta} := (E_{\widetilde H_j} -\zeta U_j)^{-1},\
\mathbf R_{j,\zeta} := P_{H}^{\widetilde H_j} R_{j,\zeta},\qquad j=1,2,\ \zeta\in \mathbb{C}: |\zeta |\not= 1. $$
From~(\ref{f2_32}),(\ref{f2_17}) it follows that
\begin{equation}
\label{f2_33}
(\mathbf R_{1,\zeta} x,y)_H = (\mathbf R_{2,\zeta} x,y)_H,\qquad x,y\in L_N,\quad
\zeta\in \mathbb{C}: |\zeta |\not= 1.
\end{equation}
Choose an arbitrary $\zeta\in \mathbb{C}: |\zeta |\not= 1$.
Since for $j=1,2$, we may write
$$ R_{j,\zeta} (E_H - \zeta A) x =
( E_{\widetilde H_j} - \zeta U_j  )^{-1} (E_{\widetilde H_j}-\zeta U_j) x = x,\qquad x\in H_0 = D(A), $$
we get
\begin{equation}
\label{f2_34}
R_{1,\zeta} u = R_{2,\zeta} u \in H,\qquad u\in H_\zeta,\ \zeta\in \mathbb{C}: |\zeta |\not= 1;
\end{equation}
\begin{equation}
\label{f2_35}
\mathbf R_{1,\zeta} u = \mathbf R_{2,\zeta} u,\qquad u\in H_\zeta,\ \zeta\in \mathbb{C}: |\zeta |\not= 1.
\end{equation}
Suppose additionally that $\zeta\not= 0$. We may write
$$ (\mathbf R_{j,\zeta} x, u)_H = (R_{j,\zeta} x, u)_{\widetilde H_j}
= (x, R_{j,\zeta}^* u)_{\widetilde H_j} $$
\begin{equation}
\label{f2_36}
= ( x, ( E_{\widetilde H_j} - R_{ j, \frac{1}{\overline{ \zeta }} } )      u)_{ \widetilde H_j} =
( x, u )_H - (x, \mathbf R_{ j, \frac{1}{\overline{ \zeta }} } u)_H,\quad
x\in L_N,\ u\in H_{ \frac{1}{ \overline{\zeta} } },\ j=1,2.
\end{equation}
Therefore we get
\begin{equation}
\label{f2_37}
(\mathbf R_{1,\zeta} x,u)_H = (\mathbf R_{2,\zeta} x,u)_H,\qquad x\in L_N,\
u\in H_{  \frac{1}{ \overline{\zeta} } },\ \zeta\in \mathbb{C}\backslash\{ 0 \}: |\zeta |\not= 1.
\end{equation}
Choose an arbitrary $\zeta\in \mathbb{C}:\ 0< |\zeta| <1$.
By~(\ref{f2_26}) an arbitrary element $y\in H$ can be represented as
$y=y_{  \frac{1}{ \overline{\zeta} }  } + y'$,
$y_{  \frac{1}{ \overline{\zeta} }  }\in H_{  \frac{1}{ \overline{\zeta} }  }$, $y'\in L_N$.
Using~(\ref{f2_33}) and~(\ref{f2_37})  we get
$$ (\mathbf R_{1,\zeta} x,y)_H = (\mathbf R_{1,\zeta} x, y_{ \frac{1}{ \overline{\zeta} } } + y')_H =
(\mathbf R_{2,\zeta} x, y_{ \frac{1}{ \overline{\zeta} } } + y')_H
= (\mathbf R_{2,z} x,y)_H, $$
where $x\in L_N,\ y\in H$.
Thus, we obtain
\begin{equation}
\label{f2_38}
\mathbf R_{1,\zeta} x = \mathbf R_{2,\zeta} x,\qquad x\in L_N,\ \zeta\in \mathbb{C}:\ 0< |\zeta| <1.
\end{equation}
Choose an arbitrary $\zeta\in \mathbb{C}:\ 0< |\zeta| <1$.
For an arbitrary $h\in H$, by~(\ref{f2_26}) we may write
$$ h = a + b,\qquad a\in L_N,\quad b\in H_\zeta. $$
Using relations~(\ref{f2_38}),(\ref{f2_35}) we obtain
$$ \mathbf R_{1,\zeta} h = \mathbf R_{1,\zeta} a + \mathbf R_{1,\zeta} b =
\mathbf R_{2,\zeta} a + \mathbf R_{2,\zeta} b = \mathbf R_{2,\zeta} h. $$
Therefore
\begin{equation}
\label{f2_39}
\mathbf R_{1,\zeta} = \mathbf R_{2,\zeta},\qquad \zeta\in \mathbb{C}:\ 0< |\zeta| <1.
\end{equation}
Observe that $\mathbf R_{1,0} = E_H = \mathbf R_{2,0}$, and the following relation
holds~\cite{cit_4000_Ch}:
$$ \mathbf R_{j,\zeta}^* = \mathbf R_{j, \frac{1}{\overline{z}} },\qquad
\zeta\in \mathbb{C}\backslash\{ 0 \}:\ |\zeta| \not= 1,\quad j=1,2. $$
Therefore
\begin{equation}
\label{f2_40}
\mathbf R_{1,\zeta} = \mathbf R_{2,\zeta},\qquad \zeta\in \mathbb{C}:\ |\zeta| \not= 1.
\end{equation}
By the inversion formula, we obtain $\mathbf{E}_{1,t} = \mathbf{E}_{2,t}$.
The obtained contradiction completes the proof.
$\Box$

We shall use the following result:
\begin{thm}~\cite[Theorem 3]{cit_4000_Ch}
\label{t2_3}
An arbitrary generalized resolvent $\mathbf{R}_\zeta$ of a closed isometric operator $U$
in a Hilbert space $H$ has the following representation:
\begin{equation}
\label{f2_41}
\mathbf R_{\zeta} =
\left[
E - \zeta ( U \oplus \Phi_\zeta )
\right]^{-1},\qquad
\zeta\in \mathbb{D}.
\end{equation}
Here $\Phi_\zeta$ is an analytic in $\mathbb{D}$ operator-valued function which values are
linear contractions (i.e. $\| \Phi_\zeta \| \leq 1$) from $H\ominus D(U)$ into $H\ominus R(U)$.

\noindent
Conversely, each analytic in $\mathbb{D}$ operator-valued function with above properties
generates by relation~(\ref{f2_41}) a generalized resolvent $\mathbf R_{\zeta}$ of $U$.
\end{thm}
Observe that relation~(\ref{f2_41}) also shows that different analytic in $\mathbb{D}$ operator-valued functions
with above properties generate different generalized resolvents of $U$.

Comparing the last two theorems we obtain the following result.
\begin{thm}
\label{t2_4}
Let the truncated matrix trigonometric moment problem~(\ref{f1_1}) be given and
condition~(\ref{f1_4}) is true. Let an operator $A$ be constructed for the
moment problem as in~(\ref{f2_11}).
All solutions of the moment problem have the following form
\begin{equation}
\label{f2_42}
M(t) = (m_{k,j}(t))_{k,j=0}^{N-1},\qquad t\in [0,2\pi],
\end{equation}
where $m_{k,j}$ are obtained from the following relation:
\begin{equation}
\label{f2_43}
\int_0^{2\pi} \frac{1}{1-\zeta e^{it}} dm_{k,j} (t) =
(\mathbf{R}_\zeta x_k,x_j)_H,\qquad z\in \mathbb{C}:\ | z | \not= 1;
\end{equation}
and
\begin{equation}
\label{f2_44}
\mathbf R_{\zeta} =
\left[
E - \zeta ( U \oplus \Phi_\zeta )
\right]^{-1},\quad
\mathbf R_{ \frac{1}{ \overline{\zeta} }} = E_H - \mathbf{R}_\zeta^*,\qquad
\zeta\in \mathbb{D}.
\end{equation}
Here $\Phi_\zeta$ is an analytic in $\mathbb{D}$ operator-valued function which values are
linear contractions from $H\ominus D(A)$ into $H\ominus R(A)$.

Conversely, each analytic in $\mathbb{D}$ operator-valued function with above properties
generates by relations~(\ref{f2_42})-(\ref{f2_44}) a solution of the moment problem~(\ref{f1_1}).

Moreover, the correspondence between all
analytic in $\mathbb{D}$ operator-valued functions with above properties
and all solutions
of the moment problem~(\ref{f1_1}) is bijective.
\end{thm}
{\bf Proof. } The proof is obvious.
$\Box$

\begin{center}
{\large\bf The truncated matrix trigonometric moment problem: the operator approach.}
\end{center}
\begin{center}
{\bf S.M. Zagorodnyuk}
\end{center}

In this paper we study the truncated matrix trigonometric moment problem.
We obtained a bijective parameterization of all solutions of this moment problem (both in nondegenerate
and degenerate cases)
via an operator approach. We use important results of
M.E.~Chumakin on generalized resolvents of isometric operators.

}
\end{document}